\documentclass[11pt]{article}
\usepackage{amsmath,amsfonts,amssymb,rotating,colordvi}
\usepackage{colordvi}
\def\qed{\nopagebreak\hfill{\rule{4pt}{7pt}}}
\def\proof{\noindent {\it{Proof.} \hskip 2pt}}

\newtheorem{theo}{Theorem}[section]

\newtheorem{conj}[theo]{Conjecture}

\numberwithin{equation}{section}

\newdimen\Squaresize \Squaresize=11pt
\newdimen\Thickness \Thickness=0.7pt
\def\Square#1{\hbox{\vrule width \Thickness
   \vbox to \Squaresize{\hrule height \Thickness\vss
    \hbox to \Squaresize{\hss#1\hss}
   \vss\hrule height\Thickness}
\unskip\vrule width \Thickness} \kern-\Thickness}

\def\Vsquare#1{\vbox{\Square{$#1$}}\kern-\Thickness}

\def\moins{\raise 1pt\hbox{{$\scriptstyle -$}}}

\begin{document}

\begin{center}
{\large \bf  Br\"{a}nd\'{e}n's Conjectures on the Boros-Moll Polynomials}
\end{center}

\begin{center}
William Y. C. Chen$^{1}$, Donna Q. J. Dou$^{2}$  and  Arthur
L. B. Yang$^{3}$\\[6pt]

$^{1,3}$Center for Combinatorics, LPMC-TJKLC\\
Nankai University, Tianjin 300071, P. R. China

$^{2}$School of Mathematics\\
Jilin University, Changchun, Jilin 130012, P. R. China

 Email: $^{1}${\tt
chen@nankai.edu.cn}, $^{2}${\tt qjdou@jlu.edu.cn}, $^{3}${\tt yang@nankai.edu.cn}
\end{center}

\noindent \emph{Abstract.}  We prove two conjectures of  Br\"{a}nd\'{e}n
on the real-rootedness of polynomials $Q_n(x)$ and $R_n(x)$ which are
related to the Boros-Moll polynomials $P_n(x)$. In fact, we show that both
$Q_n(x)$ and $R_n(x)$ form Sturm sequences.
The first conjecture  implies the $2$-log-concavity of $P_n(x)$, and the second conjecture
implies the $3$-log-concavity of $P_n(x)$.

\noindent \emph{AMS Classification 2010:} Primary 26C10; Secondary 05A20, 30C15.


\section{Introduction}

In this paper, we prove two conjectures of Br\"{a}nd\'{e}n  \cite{Bran2011}
concerning the Boros-Moll polynomials. Br\"{a}nd\'{e}n introduced
two polynomials based on the coefficients of the Boros-Moll polynomials
and conjectured that these polynomials have only real roots.
As pointed out by Br\"{a}nd\'{e}n, the first conjecture implies the $2$-fold log-concavity, or $2$-log-concavity, for short,  of the Boros-Moll polynomials, whereas the second conjecture implies the  3-log-concavity.

Let us start with some definitions. Given a finite nonnegative sequence $\{a_i\}_{i=0}^n$, we say that it is {unimodal} if there exists an integer
$m\geq 0$ such that
$$a_0\leq\cdots\leq a_{m-1}\leq
a_m\geq a_{m+1} \geq \cdots \geq a_n,$$
and we say that it is {log-concave} if
$$a_i^2-a_{i+1}a_{i-1}\geq 0$$
for $1\leq i\leq n-1$.
Define $\mathcal{L}$ to be an operator acting on the sequence $\{a_i\}_{i=0}^n$
as given by
$$\mathcal{L}(\{a_i\}_{i=0}^n)=\{b_i\}_{i=0}^n,$$
where
$b_i=a_i^2-a_{i+1}a_{i-1}$ for $0\leq i\leq n$ under the convention
that $a_{-1}=0$ and $a_{n+1}=0$.
Clearly, the sequence $\{a_i\}_{i=0}^n$ is log-concave if and only if the sequence $\{b_i\}_{i=0}^n$ is nonnegative. Given a sequence $\{a_i\}_{i=0}^n$, we say that it is
{$k$-fold log-concave}, or $k$-log-concave, if $\mathcal{L}^j(\{a_i\}_{i=0}^n)$ is a nonnegative sequence for any $1\leq j\leq k$.
A sequence $\{a_i\}_{i=0}^n$ is said to be {infinitely log-concave} if it is $k$-log-concave for all $k\geq 1$. Given a polynomial
$$f(x)=a_0+a_1x+\cdots+a_n x^n,$$
we say that $f(x)$ is log-concave (or $k$-log-concave, or infinitely log-concave) if the sequence $\{a_i\}_{i=0}^n$
of coefficients is log-concave (resp., $k$-log-concave, infinitely log-concave).

The notion of  infinite log-concavity was introduced by Boros and Moll \cite{Bormol2004} in their study of the following quartic integral
$$\int_0^{\infty}\frac{1}{(t^4+2xt^2+1)^{n+1}}{\rm d}t.$$
For any $x>-1$ and any nonnegative integer $n$, they obtained
 the following formula,
$$\int_0^{\infty}\frac{1}{(t^4+2xt^2+1)^{n+1}}{\rm d}t=\frac{\pi}{2^{n+3/2}(x+1)^{n+1/2}}P_n(x),$$
where
$$P_n(x)=\sum\limits_{j,k}{2n+1\choose 2j}{n-j\choose k}{2k+2j\choose k+j}\frac{(x+1)^j(x-1)^k}{2^{3(k+j)}}$$
are  the {Boros-Moll polynomials}.
Using Ramanujan's Master Theorem, they  derived an alternative representation of $P_n(x)$,
\begin{align}\label{eq-bm}
P_n(x)=2^{-2n}\sum_j 2^j{2n-2j\choose n-j}{n+j\choose j}(x+1)^j.
\end{align}
Write
\begin{eqnarray}\label{Boros-Moll}
P_n(x)=\sum\limits_{i=0}^n d_i(n)x^i.
\end{eqnarray}
We call  $\{d_i(n)\}_{i=0}^n$  a {Boros-Moll sequence}.
Boros and Moll proposed the following conjecture.

\begin{conj}[\cite{Bormol2004}]\label{boros-moll-conj}
The sequence $\{d_i(n)\}_{i=0}^n$ is infinitely log-concave.
\end{conj}

The log-concavity of $\{d_i(n)\}_{i=0}^n$ was  conjectured by Moll \cite{moll2002}, and it was proved by Kauers and Paule \cite{KauPau2007} by establishing  recurrence relations of the coefficients $d_i(n)$.
Chen and Xia \cite{chenxia08} showed that the polynomials $P_n(x)$
  are ratio monotone. Notice that for a positive sequence, the ratio monotone property implies both log-concavity and the spiral property.
  It is worth mentioning that there are proofs of the log-concavity without
  using recurrence relations. Llamas and Mart\'{\i}nez-Bernal \cite{llamab2010} proved that if $f(x)$ is a polynomial with nondecreasing and nonnegative coefficients, then $f(x+1)$ is log-concave. Furthermore, Chen, Yang and Zhou \cite{chyazh2010} proved that if $f(x)$ is a polynomial with nondecreasing and nonnegative coefficients, then $f(x+1)$ is ratio monotone. From (\ref{eq-bm}) it is easily seen that
    the coefficients of $P_n(x-1)$ are nondecreasing and nonnegative.
    Hence $P_n(x)$ are log-concave and ratio monotone. A combinatorial
    interpretation of the log-concavity of $P_n(x)$ has been found by
    Chen, Pang and Qu \cite{cpq2012}.

There was little progress on the higher-fold log-concavity of the Boros-Moll polynomials.
As remarked by Kauers and Paule \cite{KauPau2007}, it seems that there is little hope to prove the $2$-log-concavity of $\{d_i(n)\}_{i=0}^n$ using  recurrence relations. By constructing an intermediate function, Chen and Xia \cite{Chx2010} proved the $2$-log-concavity of $P_n(x)$ by applying recurrence relations.
Based on a technique of McNamara and Sagan \cite{mcnsag2010}, Kauers verified the infinite log-concavity of  $P_n(x)$ for $n\leq 129$.

 Br\"{a}nd\'{e}n \cite{Bran2011} presented an approach to Conjecture \ref{boros-moll-conj} by relating higher-order log-concavity to real-rooted polynomials. Boros and Moll \cite{Bormol2004} conjectured that for any nonnegative integer $n$ the sequence $\{\binom{n}{k}\}_{k=0}^n$ is   infinitely log-concave.
  Fisk \cite{fisk2008}, McNamara and Sagan \cite{mcnsag2010} and Stanley independently made the following conjecture which implies the
 conjecture of Boros and
Moll.   This conjecture has been proved by  Br\"{a}nd\'{e}n \cite{Bran2011}.

\begin{theo}  \label{bran-thm} If $f(x)=a_0+a_1x+\cdots+a_nx^n$ is a real-rooted polynomial with nonnegative coefficients,  the polynomial
$$a_0^2+(a_1^2-a_0a_2)x+\cdots+(a_{n-1}^2-a_{n-2}a_n)x^{n-1}+a_n^2x^n$$
is also real-rooted.
\end{theo}

Br\"{a}nd\'{e}n's proof is based on a symmetric function identity and the Grace-Walsh-Szeg\"o theorem concerning the location of zeros of multi-affine and symmetric polynomials. Moreover, Br\"{a}nd\'{e}n obtained a general result about the characterization of nonlinear transformations preserving real-rootedness, in the spirit of the   characterization of linear transformations preserving stability given by Borcea and Br\"{a}nd\'{e}n \cite{borbra2009}. Cardon and Nielsen \cite{carnie2011} found a combinatorial proof of Theorem \ref{bran-thm}  in terms of directed acyclic weighted planar networks. Although the Boros-Moll polynomials $P_n(x)$ are not real-rooted,  Br\"{a}nd\'{e}n \cite{Bran2011} introduced  two
polynomials related to $P_n(x)$, and conjectured that they are real-rooted.

\begin{conj}[{\cite[Conjecture 8.5]{Bran2011}}] \label{conjecture1}
For any $n\geq 1$, the polynomial
\begin{align}
Q_n(x)=\sum_{i=0}^n\frac{d_{i}(n)}{i!}x^i\label{eq-qn}
\end{align}
has only real zeros.
\end{conj}

\begin{conj}[{\cite[Conjecture 8.6]{Bran2011}}] \label{conjecture2}
For any $n\geq 1$, the polynomial
\begin{align}
R_n(x)=\sum_{i=0}^n\frac{d_{i}(n)}{(i+2)!}x^i\label{eq-rn}
\end{align}
 has only real zeros.
\end{conj}

As pointed out by Br\"{a}nd\'{e}n \cite{Bran2011},  the real-rootedness of $Q_n(x)$ implies the $2$-log-concavity of $P_n(x)$, and the real-rootedness of $R_n(x)$ implies the $3$-log-concavity of $P_n(x)$. It is worth mentioning  that Csordas \cite{csordas2011} proved the real-rootedness of some polynomials related to $Q_n(x)$. In this paper,  we shall prove the above conjectures.

\section{Proofs of Br\"{a}nd\'{e}n's Conjectures}

To prove Br\"{a}nd\'{e}n's conjectures, we shall show that
 the polynomials $Q_n(x)$ and $R_n(x)$ form
 Sturm sequences. Let us  recall a
 criterion of Liu and Wang \cite{liu-wang} which can be used to
 deduce that a polynomial sequence is a Sturm sequence.

Throughout this paper, we shall be concerned with polynomials with real coefficients. We say that a polynomial is {standard} if it is  zero or its leading coefficient is positive. Let $\mathrm{RZ}$ denote the set of polynomials with only real zeros. Suppose that $f(x)\in \mathrm{RZ}$ is a polynomial of degree $n$ with zeros $\{r_k\}_{k=1}^{n}$, and $g(x)\in \mathrm{RZ}$ is a polynomial of degree $m$ with zeros $\{s_k\}_{k=1}^{m}$.
We say that {$g(x)$ interlaces $f(x)$} if $n=m+1$ and
$$r_n\leq s_{n-1}\leq r_{n-1}\leq \cdots \leq r_2\leq s_1\leq r_1,$$
and we say that $g(x)$ strictly interlaces $f(x)$ if, in addition, they have no common zeros. We use $g(x)\preceq f(x)$ to denote that $g(x)$ interlaces $f(x)$, and use $g(x)\prec f(x)$ to denote that $g(x)$ strictly interlaces $f(x)$. For any real numbers $a, b$ and $c$, we assume that
$a\in \mathrm{RZ}$ and $a\prec bx+c$. A sequence $\{f_n(x)\}_{n\geq 0}$ of standard polynomials is said to be {a Sturm sequence} if, for $n\geq 0$,
we have $\deg f_n(x)=n$ and
$$ f_n(x)\in \mathrm{RZ} \mbox{ and } f_n(x)\prec f_{n+1}(x).$$

Liu and Wang \cite{liu-wang} gave a sufficient condition for a polynomial sequence $\{f_n(x)\}_{n\geq 0}$ to form an interlacing sequence.

\begin{theo}[{\cite[Corollary 2.4]{liu-wang}}]\label{liuwangthm2}
Let $\{f_n(x)\}_{n\geq 0}$ be a sequence of polynomials with nonnegative coefficients and $\deg f_n(x)=n$, which satisfy the following recurrence relation:
\begin{align}\label{eq-pol}
f_{n+1}(x)=a_n(x)f_{n}(x)+b_n(x)f'_{n}(x)+c_n(x)f_{n-1}(x),
\end{align}
where $a_n(x),b_n(x),c_n(x)$ are some polynomials with real coefficients.
Assume that, for some $n\geq 1$, the following conditions hold:
\begin{itemize}
\item[(i)] $f_{n-1}(x), f_n(x)\in \mathrm{RZ}$ and $f_{n-1}(x)\prec f_n(x)$; and

\item[(ii)] for any $x\leq 0$ both of $b_n(x)$ and $c_n(x)$ are nonpositive, and at least one of them is nonzero.
\end{itemize}
Then we have $f_{n+1}(x)\in \mathrm{RZ}$ and $f_n(x)\prec f_{n+1}(x)$.
\end{theo}

To prove Conjectures \ref{conjecture1} and \ref{conjecture2},
 we proceed to derive recurrence relations for $Q_n(x)$ and $R_n(x)$ based on the recurrence relations of the coefficients $d_i(n)$ of the Boros-Moll
polynomials $P_n(x)$.  Kauers and Paule \cite{KauPau2007} proved that
\begin{align}
d_i(n+1)& =\frac{n+i}{n+1}d_{i-1}(n)+\frac{4n+2i+3}{2(n+1)}d_{i}(n),\quad 0\leq i\leq n+1, \label{rec1}\\[6pt]
d_i(n+2)& =\frac{8n^2+24n+19-4i^2}{2(n+2-i)(n+2)}d_{i}(n+1)\nonumber\\
        &\qquad -\frac{(n+i+1)(4n+3)(4n+5)}{4(n+2-i)(n+1)(n+2)}d_{i}(n), \quad 0\leq i\leq n+1.\label{rec2}
\end{align}
In fact, \eqref{rec1} can be easily derived from \eqref{rec2}.
Note that
Moll \cite{moll2007} independently derived the relation \eqref{rec2} via the WZ-method.

\begin{theo}\label{rec-qm}
For $n\geq 1$, we have the following recurrence relation
\begin{align}
Q_{n+1}(x)=&\left(\frac{(2n+1)x}{(n+1)^2}+\frac{8n^2+8n+3}{2(n+1)^2}\right)Q_n(x)\nonumber\\[6pt]
&-\frac{(4n-1)(4n+1)}{4(n+1)^2}Q_{n-1}(x)+\frac{x}{(n+1)^2}Q'_n(x). \label{eq-rec-qm}
\end{align}
\end{theo}

\proof For $n\geq 1$,  relation \eqref{eq-rec-qm} can be rewritten as 
\begin{align}
4(n+1)^2d_i(n+1)&=2(8n^2+8n+3+2i)d_{i}(n)+{4i(2n+1)d_{i-1}(n)} \nonumber\\
& \qquad -(16n^2-1)d_{i}(n-1), \label{eq-inter-q}
\end{align}
where $0\leq i\leq n+1$.
From \eqref{rec1} it follows  that
\begin{align}\label{eq-new}
d_{i-1}(n)=\frac{n+1}{n+i}d_i(n+1)-\frac{4n+2i+3}{2(n+i)}d_{i}(n).
\end{align}
Substituting \eqref{eq-new} into \eqref{eq-inter-q}, we get
\begin{align}
d_i(n+1)& =\frac{8n^2+8n+3-4i^2}{2(n+1-i)(n+1)}d_{i}(n)\nonumber\\
        &\qquad -\frac{(n+i)(4n-1)(4n+1)}{4n(n+1)(n+1-i)}d_{i}(n-1).\label{eq-inter}
\end{align}
It is easily checked that the above relation (\ref{eq-inter})
coincides with
\eqref{rec2} with $n$ replaced by $n-1$.
This completes the proof.
\qed

Using the above recurrence relation and the criterion of Liu and Wang,
 we can deduce that the polynomials $Q_n(x)$ form a Sturm sequence.
 This leads to an affirmative answer to Conjecture \ref{conjecture1}.

\begin{theo} \label{thm-1} The polynomial sequence $\{Q_n(x)\}_{n\geq 0}$ is a Sturm sequence.
\end{theo}

\proof Clearly, we have $\deg (Q_n(x))=n$.
It suffices to prove that $Q_n(x)\in \mathrm{RZ}$ and $Q_n(x)\prec Q_{n+1}(x)$ for any $n\geq 0$. We use induction on $n$.
By convention,
$$Q_0(x),Q_1(x)\in \mathrm{RZ} \quad \mbox{and}\quad Q_0(x)\prec Q_1(x).$$
Assume that
$$Q_{n-1}(x),Q_n(x)\in \mathrm{RZ} \quad \mbox{and}\quad Q_{n-1}(x)\prec Q_n(x).$$
We proceed to verify that
$$Q_{n+1}(x)\in \mathrm{RZ} \quad \mbox{and}\quad Q_n(x)\prec Q_{n+1}(x).$$
 We see that the recurrence relation \eqref{eq-rec-qm} of $Q_n(x)$ is of the form \eqref{eq-pol}  in  Theorem \ref{liuwangthm2}, where the polynomials $a_n(x),b_n(x),c_n(x)$  are given by
\begin{align*}
a_n(x)&=\frac{(2n+1)x}{(n+1)^2}+\frac{8n^2+8n+3}{2(n+1)^2},\\
b_n(x)&=\frac{x}{(n+1)^2},\\
c_n(x)&=-\frac{(4n-1)(4n+1)}{4(n+1)^2}.
\end{align*}
For $n\geq 1$ and $x\leq 0$, one can check that
$$b_n(x)\leq 0  \quad \mbox{and}\quad  c_n(x)<0.$$
In view of Theorem \ref{liuwangthm2}, we find that $Q_{n+1}(x)\in \mathrm{RZ}$ and $Q_n(x)\prec Q_{n+1}(x)$. This completes the proof. \qed

The following recurrence  relation for $R_n(x)$ can be proved in a way similar to the proof of Theorem \ref{rec-qm}.

\begin{theo}\label{rec-rm}
For $n\geq 1$, we have 
\begin{align}
R_{n+1}(x)&=\left(\frac{(2n+1)x}{(n+1)(n+3)}+\frac{8n^2+8n+7}{2(n+1)(n+3)}\right)R_n(x)\nonumber\\[6pt]
&\qquad -\frac{(4n-1)(4n+1)(n-2)}{4n(n+1)(n+3)}R_{n-1}(x)+\frac{5x}{(n+1)(n+3)}R'_n(x).\label{eq-rec-rm}
\end{align}
\end{theo}

Using the above recurrence relation, we obtain the following theorem, 
which leads to an affirmative answer to Conjecture \ref{conjecture2}.

\begin{theo} \label{thm-2} The polynomial sequence $\{R_n(x)\}_{n\geq 0}$ is a Sturm sequence.
\end{theo}

\proof The proof is analogous to that of Theorem  \ref{thm-1}.
It is routine to verify that
$$R_0(x),R_1(x),R_2(x),R_3(x)\in \mathrm{RZ} \quad \mbox{and}\quad R_0(x)\prec R_1(x)\prec R_2(x)\prec R_3(x).$$
It remains to show that $R_n(x)\in \mathrm{RZ}$ and $R_{n-1}(x)\prec R_{n}(x)$ for $n\geq 3$. We use induction $n$.
Assume that
$$R_{n-1}(x),R_n(x)\in \mathrm{RZ} \quad \mbox{and}\quad R_{n-1}(x)\prec R_n(x).$$
We wish to prove that
$$R_{n+1}(x)\in \mathrm{RZ} \quad \mbox{and}\quad R_n(x)\prec R_{n+1}(x).$$
The recurrence relation \eqref{eq-rec-rm} of $R_n(x)$ is  of the form \eqref{eq-pol}  in  Theorem \ref{liuwangthm2}, and the polynomials $a_n(x),b_n(x),c_n(x)$  are given by
\begin{align*}
a_n(x)&=\frac{(2n+1)x}{(n+1)(n+3)}+\frac{8n^2+8n+7}{2(n+1)(n+3)},\\
b_n(x)&=\frac{5x}{(n+1)(n+3)},\\
c_n(x)&=-\frac{(4n-1)(4n+1)(n-2)}{4n(n+1)(n+3)}.
\end{align*}
For $n\geq 3$ and $x\leq 0$, we find that
$$b_n(x)\leq 0  \quad \mbox{and}\quad  c_n(x)<0.$$
By Theorem \ref{liuwangthm2}, we conclude that $R_{n+1}(x)\in \mathrm{RZ}$ and $R_n(x)\prec R_{n+1}(x)$. This completes the proof. \qed

\vskip 3mm
\noindent {\bf Acknowledgments.} This work was supported by the 973 Project, the PCSIRT Project
of the Ministry of Education, and the National Science Foundation of China.

\end{document}